\documentclass[a4paper]{article}

\usepackage[english]{babel}
\usepackage[utf8]{inputenc}
\usepackage{amsmath}
\usepackage{amsfonts}
\usepackage{amsthm}
\usepackage{graphicx}
\usepackage[colorinlistoftodos]{todonotes}
\usepackage{amssymb}

\usepackage{xy}
\input xy
\xyoption{all}

\def\acts{\curvearrowright}

\newtheorem{theorem}{Theorem}[section]
\newtheorem{lemma}[theorem]{Lemma}
\newtheorem{proposition}[theorem]{Proposition}
\newtheorem{corollary}[theorem]{Corollary}

\newtheorem{example}[theorem]{Example}

\theoremstyle{definition}
\newtheorem{definition}[theorem]{Definition}
\theoremstyle{remark}

\title{Warped cones and proper affine isometric actions of discrete groups on Banach spaces
			\thanks{Supported by NSFC, Nos. 11231002, 11420101001.}}

\author{Qin Wang and Zhen Wang}

\date{\today}

\begin{document}
\maketitle

\begin{abstract}
Warped cones are metric spaces introduced by John Roe from discrete group actions on compact metric spaces to produce interesting examples in coarse geometry.
We show that a certain class of warped cones $\mathcal{O}_\Gamma (M)$ admit a fibred coarse embedding into a $L_p$-space ($1\leq p<\infty$) if and only if the discrete group $\Gamma$
admits a proper affine isometric action on a $L_p$-space. This actually holds for any class of Banach spaces stable under taking Lebesgue-Bochner $L_p$-spaces and ultraproducts, e.g.,
uniformly convex Banach spaces or Banach spaces with nontrivial type. It follows that the maximal coarse Baum-Connes conjecture and the coarse Novikov conjecture hold for a certain class of
warped cones which do not coarsely embed into any $L_p$-space for any $1\leq p<\infty$.
\end{abstract}

\section{Introduction}

Warped cones were introduced by John Roe \cite{Roe95, Roe05}. Typically, a warped cone $\mathcal{O}_\Gamma(M)$ is constructed from the action of a finitely generated group $\Gamma$ on a compact Riemannian
manifold $M$ by diffeomorphisms. It is a metric space whose underlying topological space is the open cone $\mathcal{O} M$ on $M$, but whose coarse geometry produces large groups of translations. It was shown in
\cite{Roe05} that the warped cone construction is sufficiently flexible to create examples of metric spaces which do or do not have Yu's property A, or which do not admit a coarse embedding into
Hilbert space. Coarse geometry of warped cones over rotations on the circle was studied afterwards in \cite{Kim}.
\par
Recently, warped cones were further studied in several very interesting works \cite{DruNow,NowSaw,Saw,Vigolo,Saw17}. In \cite{DruNow}, C. Dru\c{t}u and P. Nowak constructed ghost projections for warped cones
from certain actions with a spectral gap, and conjectured that such warped cones provide new counterexamples to the coarse Baum-Connes conjecture.  In \cite{Saw}, D. Sawicki provided
first examples of warped cones which coarsely embed into Hilbert space but do not have property A and, along the way, founded out sevaral formulas simplifying calculation of the warped metric.  In \cite{NowSaw}, P. Nowak and D. Sawicki showed that, due to spectral gaps of group actions, warped cones like
$\mathcal{O}_{\mathbb{F}_2}(\text{SU(2)})$ does not coarsely embed into any $L_p$-space for any $1\leq p<\infty$. In \cite{Vigolo}, F. Vigolo discovered that a large class of
warped cones in which the group action is expanding in measure contain a family of expanders and hence do not coarsely embed into any $L_p$-space. In \cite{Saw17}, D. Sawicki showed that for many sources of actions $\Gamma\acts Y$ and any Banach space $X$, 
the levels of the warped cone $\mathcal{O}_{\Gamma}(Y)$ are expanders with respect to $X$ if and only if the induced representation of $\Gamma$ on $L_2(Y; X)$ has a spectral gap.
\par
On the other hand, a notion of {\em fibred coarse embeddings} for metric spaces was introduced by X. Chen, Q. Wang and G. Yu in \cite{CWY13} as a far-reaching generalization of Gromov's notion of coarse embeddings. It was
shown in \cite{CWY13} that the maximal coarse Baum-Connes conjecture holds for metric spaces of bounded geometry which admit a fibred coarse embedding into Hilbert space. This class of spaces includes all
expander graphs with large girth and all box spaces of residually finite groups with the Haagerup property. It was shown in \cite{CWW} that a box space $\Box_{\{\Gamma_n\}} \Gamma$ of a residually finite group $\Gamma$
admits a fibred coarse embedding into Hilbert space if and only if the group $\Gamma$ has the Haagerup property. Furthermore, it was actually shown in \cite{CWW}  that if $\Gamma$ admits
a proper isometric action on a metric space $Y$, then any box space  $\Box_{\{\Gamma_n\}} \Gamma$  admits a fibred coarse embadding into $Y$. Very recently, S. Arnt showed in \cite{Arnt} that the converse of
this last statement for $L_p$ spaces (or other stable classes of Banach spaces) is also true. Consequently, a box space $\Box_{\{\Gamma_n\}} \Gamma$  of a residually finite group $\Gamma$ admits a fibred coarse embedding
in a $L_p$ space ($1\leq p<\infty$) if and only if the group $\Gamma$ admits a proper affine isometric action on a $L_p$ space.
\par
In this paper, we shall investigate the connection between fibred coarse embeddability of warped cones and existence of proper affine isometric actions of discrete groups on Banach spaces. Note that a large class of transformation groups $(\Gamma, M)$
admits a linearization by unitary operators in a Hilbert space, see Section 2.1 below for details. This allows us to show in Section 3 that, for such type of dynamics, if $\Gamma$ admits a proper affine isometric action on a $L_p$-space or a $\ell_p$-space, then the warped cone  $\mathcal{O}_\Gamma(M)$ admits a fibred coarse embedding into a $L_p$-space or a $\ell_p$-space, respectively. It follows that, for example, the warped cone  $\mathcal{O}_{\mathbb{F}_2}(\text{SU(2)})$  admits a fibred coarse embedding into Hilbert space and, by applying in addition with G. Yu's results in \cite{Yu05}, that
certain warped cones of actions by hyperbolic groups admit a fibred coarse embedding into an $\ell_p$ space with large $p>2$. In sections 4, we discuss the converse direction of the above results, namely, the implication from fibred coarse embeddability of warped cones to proper affine isometric actions of discrete groups on Banach spaces.
The results in Section 3 and Section 4 work very well in the case of finitely generated dense subgroups $\Gamma$ of compact metrizable groups $G$. These are specified in Section 5 where we summarize that, for example,
the warped cone $\mathcal{O}_{\Gamma}(G)$ admits a fibred coarse embedding into a $L_p$-space if and only if $\Gamma$ admits a proper affine isometric action on a $L_p$-space. The class of $L_p$ spaces here can be replaced by any other class of Banach spaces $E$ which are stable under taking ultraproducts and taking Lebesgue-Bochner spaces $L_p(\Omega, \mu; E)$. In section 6, we apply results of previous sections to higher index problems on warped cones. It turns out that the maximal coarse Baum-Connes conjecture holds for a large class of warped cones, including e.g. $\mathcal{O}_{\mathbb{F}_2}(\text{SU(2)})$, which do not coarsely embed any $L_p$ space for any $1\leq p<\infty$, and the coarse Novikov conjecture holds for a larger class of warped cones, including warped cones over certain actions by hyperbolic groups.

\section{Warped cones and linearization of transformation groups}

% \subsection{Warped metric}
In this section, we shall recall the notions of warped metric and warped cones introduced by John Roe \cite{Roe05}, together with some results in literature on linearization of transformation groups which we need in construction of warped cones for group actions on general compact metrizable spaces.
\par
% \begin{definition}[\cite{Roe05}]
Let $(X, d)$ be a proper metric space and let $\Gamma$ be a group acting by homeomorphisms on $X$, provided with a finite generating set $S$. The {\em warped metric} $\delta_\Gamma$ on $X$ \cite{Roe05} is the greatest metric that satisfies the inequalities
$$\delta_\Gamma(x, x')\leq d(x, x'), \quad \delta_\Gamma(x, sx)\leq 1 \;\; \forall s\in S. $$
%\end{definition}
%\begin{proposition}
For $\gamma\in \Gamma$, let $|\gamma|$ denote the word length of $\gamma$ relative to the generating set $S$. Let $x, y\in X$. The warped distance from $x$ to $y$ is the infimum of all sums \cite{Roe05}
$$\sum d(\gamma_ix_i, x_{i+1})+|\gamma_i|$$
taken over all finite sequences $x=x_0, x_1, \cdots, x_n=y$ in $X$ and $\gamma_0, \gamma_1, \cdots, \gamma_{n-1}$ in $\Gamma$.
%\end{proposition}
It was shown in \cite{Roe05} that the warped metric is a proper metric. The coarse structure induced by the warped metric does not depend on the choice of the generating set $S$ for $\Gamma$, nor on the choice of metric $d$ within the coarse structure of $X$.

% \subsection{Warped cones on compact manifolds}

% \begin{definition}[\cite{Roe05} Open cones on smooth compact manifolds]
Let $M$ be a smooth compact manifold. The {\em open cone $X=\mathcal{O} M$} on $M$ is defined in  two equivalent ways \cite{Roe05}:
\begin{description}
\item[(a)] Embed $M$ smoothly in a high-dimensional sphere $S^{N-1}$, and let $X$ be the union of all the rays through the origin of $\mathbb{R}^N$ that meet the embedded copy of $M$; equip $X$ with the metric (distance function) induced from $\mathbb{R}^N$.
\item[(b)] Let $X=M\times [1, \infty)$ as a manifold, equipped with the Riemannian metric $t^2g_M+g_{\mathbb{R}}$, where $g_{M}$ is the Riemannian metric on $M$, and $g_{\mathbb{R}}$ is the standard Riemannian metric on $\mathbb{R}$, and $t$ is the coordinate on $\mathbb{R}$.
\end{description}
% \end{definition}
\par
% \begin{remark}
It is straightforward to show that up to coarse equivalence these constructions are independent of the choices involved and, moreover, that both constructions yield coarsely equivalent spaces $X$.
% \end{remark}
\par
\begin{definition} [Warped cones on smooth compact manifolds \cite{Roe05}]
Let $M$ be a smooth compact manifold, and let $\Gamma$ be a finitely generated group acting on $M$ by homeomorphisms. The {\em warped cone} $\mathcal{O}_\Gamma (M)$ is the coarse space obtained by warping the open cone $\mathcal{O}M$ along the induced $\Gamma$-actions.
\end{definition}
\par
% \begin{remark} [\cite{Roe05}]
% Suppose that $\Gamma$ acts on $M$ by Lipschitz homeomorphisms. Then the warped cone $\mathcal{O}_\Gamma (M)$ has bounded geometry.
% \end{remark}
\par

% \subsection{Warped cones on general compact metric spaces}
\par
For group actions on general compact metric spaces $Y$, there are two analogous ways to construct open cones and warped cones.
\par

% \begin{definition} [\cite{Saw} Intrinsic open cones]
Let $(Y, d_Y)$ be a compact metric space. The {\em (intrinsic) open cone of $Y$} \cite{DruNow,Saw}, also denoted by $\mathcal{O} Y$, is the space $Y\times [1, +\infty)$ with the metric $d_{\mathcal{O} Y}$ defined by
$$d_{\mathcal{O} Y} \big( (y, t), (y', t') \big)=|t-t'|+\big( \min{(t, t') \big) \cdot d_Y (y, y')}.$$
% \end{definition}
\par
\begin{definition}  [Warped cones on general compact metric spaces \cite{DruNow,Saw}]
\label{warped cones on general}
Let $(Y, d_Y)$ be a compact metric space and let $\Gamma$ be a finitely generated group acting on $Y$ by homeomorphisms. The {\em warped cone of $Y$}, denoted by $\mathcal{O}_\Gamma (Y)$, is the (intrinsic) open cone
$\mathcal{O} Y$ with the warped metric, where the warping group action is defined by $\gamma \cdot (y, t) = (\gamma y, t)$.
\end{definition}
\par
D. Sawicki \cite{Saw} provided the following formulas:
\par
\begin{lemma}[\cite{Saw}]
For $t\geq 1$, let $\delta_t$ denote the warped metric on $(Y, \;  t d_Y)$. Then the warped metric $\delta_\Gamma$ on the cone $\mathcal{O} Y$ is given by the formula
$$\delta_\Gamma \Big( (y, t), (y', t+s) \Big) = s + \delta_t (y, y')$$
If $\Gamma$ acts on $(X, d_X) $ by isometries, then the warped metric can be vastly simplified:
$$\delta_\Gamma (x, x')=\inf_{\gamma\in \Gamma} \Big( |\gamma|+ d_X (\gamma x, x') \Big).$$
\end{lemma}
\par
Another way to construct open cones is to embed the compact metric space into a separable Hilbert space, which is always possible.
\par
\begin{definition} [Open cones via embedding in Hilbert spaces \cite{Will09}]
Let $(Y, d_Y)$ be a compact metric space. Embed $Y$ into the unit sphere of a separable Hilbert space $H$. Then the {\em open cone $\mathcal{O} Y$ of $Y$ in $H$} is the union of all rays through the origin and some point of the embedded copy of $Y$, equipped with the induced metric (distance function) from the Hilbert space $H$.
\end{definition}
\par
Note that distinct embeddings can give rise to coarsely distinct open cones. For purpose of this paper, we need stronger embeddings taking into account of group actions, as specified in the following subsection.

\subsection{Linearization of transformation groups}
Let $G$ be a topological transformation group acting on a metrizable space $M$. It is well known that $M$ can be topologically embedded in a Euclidean space $\mathbb{R}^N$ if it is finite dimensional,  and in a Hilbert space $H$ in all cases. What happens to the action of $G$ in the process? The following notion was introduced in \cite{BaayenGroot}:
\par
\begin{definition} [P. C. Baayen and J. De Groot \cite{BaayenGroot}]
Let $G$ be a topological transformation group acting on a metrizable space $M$. If $M$ can be embedded in a Hilbert space $H$ in such a way that the transformations $\gamma\in G$ become restrictions of continuous linear mappings of $H$ onto it self, we say that the system {\em $(G, M)$ can be linearized in $H$}.
\end{definition}
\par
There were extensive researches in history in improving the behavior of transformations. In the following, we only cite several results relevant to this paper. In 1957, G. D. Mostow proved the following theorem.
\par
\begin{theorem} [G. D. Mostow \cite{Mostow1957}]
If $M$ is a compact manifold and $G$ is a compact Lie group of homeomorphisms of $M$, then the action of $G$ can be linearized by orthogonal transformations of a Euclidean space $\mathbb{R}^N$.
\end{theorem}
\par
Recall {\em the Nash embedding theorem}  \cite{Nash} asserts that any Riemannian manifold possesses an isometric embedding into a Euclidean space of sufficiently large dimension
(isometric means preserving the length of every path).
In 1980, J. D. Moore and R. Schlafly proved an equivariant version of Nash's embedding theorem as follows:
\par
\begin{theorem}[J. D. Moore and R. Schlafly \cite{MooreSchlafly}]
If $M$ is a compact Riemannian manifold and $G$ is a compact Lie group which acts on $M$ by isometries, there is an orthogonal representation $\rho$ of $G$ on some Euclidean space $\mathbb{R}^N$ and an isometric embedding from $M$ into $\mathbb{R}^N$ which is equivariant with respect to $\rho$.
\end{theorem}
\par
In the above theorem, the representation $\rho$ can be regarded as a Lie group homomorphism from $G$ into the orthogonal group $O(N)$ which acts on $\mathbb{R}^N$ by rotations and reflections. A smooth map $f:M\to \mathbb{R}^N$ is equivariant with respect to $\rho$ if and only if $f(\gamma x)=\rho(\gamma) f(x)$, for all $\gamma\in G$ and $x\in M$.
\par
For more general transformation groups, P. C. Baayen and J. de Groot showed in 1968 that, for an extensive class of locally compact transformation groups $(G, M)$, there exit a topological embedding of $M$ in a Hilbert space $H$ such that the transformations $\gamma\in G$ become restrictions of invertible bounded linear operations in $H$. In particular, if $G$ is compact, it is even possible to linearize $G$ by unitary operators of $H$.
\par
\begin{theorem}[P. C. Baayen and J. de Groot \cite{BaayenGroot}]
Let $G$ be a topological transformation group acting on a metrizable space $M$. If $G$ is either compact, or countable, or locally compact abelian, or compactly generated, its action can be linearized in a Hilbert space $H$. If $G$ is a compact, it is possible to liearize $G$ by unitary transformations in a Hilbert space $H$.
\end{theorem}
\par
%\subsection{Warped cones in linearizations by unitary operators}
In this paper, we shall define warped cones over group actions on general {\em compact metrizable spaces} by taking into account of linearization of the transformation groups. Let $M$ be a compact metrizable space, and let $\Gamma$ be a finitely generated group acting on $M$ by homeomorphisms. Suppose $(\Gamma, M)$ can be linearized by unitary operators in a Hilbert space $H$. By identifying $H$ with a hyperplane of $H\oplus \mathbb{C}$ by the mapping $\xi\mapsto (\xi, 1)$ if necessary, we may assume that $M$ is embedded into the unit sphere of $H$. Then the open cone $\mathcal{O}M $ of $M$ is the subspace of $H$ consisting of all the rays through the origin that meet $M$, and the group $\Gamma$ acts on $\mathcal{O} M$ by restrictions of unitary operators in $H$.
\par
\begin{definition}[Warped cones in linearizations by unitary operators]
Let $M$ be a compact metrizable space, and let $\Gamma$ be a finitely generated group acting on $M$ by homeomorphisms. Suppose $(\Gamma, M)$ can be linearized by unitary operators in a Hilbert space $H$.
The {\em warped cone of $(\Gamma, M)$ in the linearization by unitary operators in $H$} is defined to be the open cone $\mathcal{O} M$ in $H$ with the warped metric $\delta_\Gamma$ along the action of $\Gamma$ on $\mathcal{O} M$ obtained by restrictions of the unitary representation of $\Gamma$ on $H$.
\end{definition}
\par
\begin{example}
Let $\Gamma$ be a residually finite group, and let $\{\Gamma_n\}^\infty_{n=1}$ be a nested sequence of finite index normal subgroups with trivial intersection. Let $G:=\widehat{\Gamma}((\Gamma_n))$ be the corresponding profinite completion. Then $G$ is a compact metrizable group containing $\Gamma$ as a dense subgroup. It follows that the system $(\Gamma, G)$ can be linearized in a Hilbert space $H$ by unitary operators. The warped cone $\mathcal{O}_\Gamma (\widehat{\Gamma}((\Gamma_n)))$ in linearization by unitary operators is defined as above.
\end{example}

\section{Fibred coarse embeddability of warped cones into Banach spaces}
In this section, we shall study fibred coarse embeddability into Banach spaces for a class of warped cones in which the acting group admits a proper affine isometric action on a Banach space. We first recall the concept of {\em fibred coarse embedding} into Banach space for a metric space, introduced in \cite{CWY13}, which is a generalization of Gromov's notion of
coarse embedding \cite{Grom93}. Let $E$ be a Banach space, as a model space.
\par
\begin{definition}[\cite{CWY13}]
\label{fibred}
A metric space $(X, d)$ is said to admit a {\em fibred coarse embedding into $E$} if there exist
\begin{itemize}
\item a field of Banach spaces $(E_x)_{x\in X}$ over $X$;
\item a section $s: X\to \bigsqcup_{x\in X} E_x$ (i.e. $s(x)\in E_x$);
\item two non-decreasing functions $\rho_1$ and $\rho_2$ from $[0, \infty)$ to $(-\infty, \infty)$ with $\lim_{r\to \infty}\rho_i(r)=\infty$ ($i=1, 2$)
\end{itemize}
such that for any $r>0$ there exists a bounded subset $K\subset X$ for which there exists a ``trivialization''
$$t_C: (E_x)_{x\in C}\longrightarrow C\times E$$
for each subset $C\subset X\backslash K $ of diameter less than $r$, i.e. a map from $(E_x)_{x\in C}$ to the constant field $C\times E$ over $C$ such that the restriction  of $t_C$ to
the fiber $E_x$ ($x\in C$) is an affine isometry $t_C(x): E_x\to E$, satisfying
\begin{itemize}
\item[(1)] for any $x, y\in C$, $\rho_1(d(x, y))\leq \|t_C(x)(s(x))-t_C(y)(s(y))\|\leq \rho_2(d(x, y))$;
\item[(2)] for any two subsets $C_1, C_2\subset X\backslash K$ of diameter less than $r$ with $C_1\cap C_2\not=\emptyset$, there exists an affine isometry $t_{C_1C_2}: E\to E$ such that
		$t_{C_1}(x)\circ t_{C_2}^{-1}(x)=t_{C_1C_2}$ for all $x\in C_1\cap C_2$. \hfill{$\heartsuit$}
\end{itemize}
\end{definition}
\par
Before stating the main result of this section, we fix some notations. Let $E$ be a Banach space and $H$ a Hilbert space. For any $1\leq p<\infty$, we denote by $E\oplus_p H$ the $\ell_p$-direct sum:
$$E\oplus_p H=\{(\xi, \eta) | \xi\in E, \eta\in H\}$$
with the norm
$$\|(\xi, \eta)\|_p = \big( \|\xi\|^p_E+\|\eta\|_H^p \big)^{\frac{1}{p}}.$$
Let $\Gamma$ and $X$ be metric spaces. We denote by $\Gamma\times_1 X$ the $\ell_1$-product space: $\Gamma\times_1 X=\{(\gamma, x)| \gamma\in \Gamma, x\in X\}$ with the metric
$$d_1((\gamma_1, x_1), (\gamma_2, x_2))=d(\gamma_1, \gamma_2)+d(x_1, x_2).$$
Note that if $\Gamma$ and $X$ admit coarse embedding into $E$ and $H$, respectively, then $\Gamma\times_1 X$ admits a coarse embedding into $E\oplus_p H$. (A more general situation for coarse embedding product metric spaces into
$\ell_p$-spaces was carefully presented by D. Sawicki in Proposition A.3 in \cite{Saw}.)
\par
The main result of this section is the following theorem.
\par
\begin{theorem}
\label{group to cone}
Let $M$ be a compact metrizable space and let $\Gamma$ be a finitely generated group acting on $M$ by homeomorphisms. Suppose
\begin{enumerate}
\item[(1)] the system $(\Gamma, M)$ admits a linearization in a Hilbert space $H$ by unitary operators;
\item[(2)] $M$ contains a dense and free orbit of $\Gamma$;
\item[(3)] $\Gamma$ admits a proper affine isometric action on a Banach space $E$.
\end{enumerate}
Then, for any $1\leq p < + \infty$, the warped cone $\mathcal{O}_\Gamma (M)$ admits a fibred coarse embedding into the Banach space $E\oplus_p H$.
\end{theorem}
\par
{\bf Proof.} Since the system $(\Gamma, M)$ can be linearized in the Hilbert space $H$ by unitary operators, we may assume that the open cone $\mathcal{O} M$ is a subset of $H$,  and $\gamma\in \Gamma$ acts on $\mathcal{O} M$ by restrictions of unitary operators $U_\gamma$ on $H$, where $\gamma\mapsto U_\gamma$ denotes the unitary representation of $\Gamma$ on $H$ in the linearization.
\par
Let $\alpha$ be the proper affine isometric action of $\Gamma$ on the Banach space $E$. Then there exist a linear isometric representation $\pi$ of $\Gamma$ on $E$ and a $1$-cocycle $b:\Gamma\to E$ for $\pi$ such that
$$\alpha_\gamma (\xi) = \pi_\gamma (\xi) + b_\gamma$$
for all $\gamma\in \Gamma$ and $\xi\in E$. The action $\alpha$ is metrically proper in the sense $\lim_{|\gamma|\to \infty} \|b_\gamma\| =+\infty.$
\par
If follows that the map
\[
\begin{array}{ccc}
\Gamma \times_1 \mathcal{O} M & \longrightarrow & E \oplus_p H  \\
(\gamma, x) & \longmapsto & (b_{\gamma^{-1}}, x)
\end{array}
\]
is a coarse embedding. Namely, there exist non-decreasing maps $\rho_1$ and $\rho_2$ on $[0, +\infty)$ tending to infinity,  such that for any pairs of elements
$(\gamma_1, x_1), (\gamma_2, x_2)\in \Gamma\times_1 \mathcal{O} M$, we have
$$\big\| ( b_{\gamma_1^{-1}}, x_1) - ( b_{\gamma_2^{-1}}, x_2) \big\|_p \in [\rho_1, \rho_2] \Big( d_{1} \big( (\gamma_1, x_1), (\gamma_2, x_2) \big) \Big),$$
i.e.,
$$
 \Big( \| b_{\gamma_1^{-1}} - b_{\gamma_2^{-1}} \|_E ^p + \| x_1- x_2 \|_H ^p \Big)^{\frac{1}{p}} \in [\rho_1, \rho_2] \Big( | \gamma_1 \gamma_2^{-1}| +\|x_1-x_2\|_H \Big).   \eqno{(*)}
$$
\par
To show that the warped cone $\mathcal{O}_\Gamma (M)$ admits a fibred coarse embedding into $E\oplus_p H$, we take the control functions $\rho_1$ and $\rho_2$ as above, and a field of Banach spaces to be
$(E\oplus_p H)_{x\in \mathcal{O}_\Gamma (M)}$, where $(E\oplus_p H)_x=E\oplus_p H$ for all $x\in \mathcal{O}_\Gamma (M)$.
\par
A section
$$s: \mathcal{O}_\Gamma (M) \to \bigsqcup_{x\in \mathcal{O}_\Gamma(M)} ( E\oplus_p H)_{x} $$
is defined as follows.
\par
To start with, since $M$ contains a dense and free orbit of $\Gamma$ by condition (2), there exists a point $P_0\in M$ such that the orbit $\Gamma P_0$ is dense in $M$,  and $\gamma_1 P_0=\gamma_2 P_0$ implies $\gamma_1=\gamma_2$ for any $\gamma_1, \gamma_2\in \Gamma$. To simply notation, without loss of generality, in the following we shall abuse notion to {\bf assume
$$M=\Gamma P_0, $$
}
so that for any $x\in M$ there exists a unique $\gamma_x\in \Gamma$ such that $x=\gamma_x P_0$.
\par
For $t>1$, denote by $\mathcal{O}_\Gamma (M) _t$ and $\mathcal{O}_\Gamma (M) _{[t, \infty)}$ the $t$-level set and the $[t, \infty)$-part of the warped cone, respectively. Denote $P_0^t:=(P_0, t)$ in $\mathcal{O}_\Gamma (M) _t$. Then for any
$x\in \mathcal{O}_\Gamma (M) _t$, there exists a unique $\gamma_x\in \Gamma$ as above such that $x=\gamma_x \cdot P_0^t$. Now we define the section
$$s: \mathcal{O}_\Gamma (M) \longrightarrow \bigsqcup_{x\in \mathcal{O}_\Gamma(M)} ( E\oplus_p H)_{x} $$
by the formula:
$$x=\gamma_x \cdot P_0^{t(x)} \longmapsto (b_{\gamma_x^{-1}}, x).$$
\par
Next, we construct ``local trivializations'' outside bounded subsets of the warped cone. For any $R>0$, choose $t_R>0$ large enough such that for any $x\in \mathcal{O}_\Gamma (M)_{[t_R, \infty)}$ we have
$$ \text{Ball}_{\mathcal{O}_\Gamma (M)} \big( x, 10R \big)  \; \subset \bigsqcup_{|g|\leq 10R} \text{Ball}_{\mathcal{O}M} \big( gx, 10R \big),$$
and the balls $\text{Ball}_{\mathcal{O} M} \big( gx, 10R \big)$ in the open cone ${\mathcal{O} M}$ for the group elements $g\in \Gamma$ with $|g|\leq 10R $ are mutually disjoint.
\par
Let $K={\mathcal{O}_\Gamma (M)}_{[1, t_R)}$.   For any $C\subset X\backslash K= \mathcal{O}_\Gamma (M)_{[t_R, \infty)} $ of diameter less than $R$, choose and fix a point $z_C\in \mathcal{O}_\Gamma(M)_{[t_R, \infty)}$ such that
$$C\subset \text{Ball}_{\mathcal{O}_\Gamma (M)} (z_C, 3R) \subset \bigsqcup _{|g|\leq 3R} \text{Ball}_{\mathcal{O}(M)} (g z_C, 3R).$$
Then, for any $x\in C$, there exist a unique $z\in \text{Ball}_{\mathcal{O} M} (z_C, 3R)$ and a unique $g\in \Gamma$ such that $x=gz$. Note that $x$ and $z$ are in the same $t$-level set of the warped cone, i.e. $t(x)=t(z)$, where,  for a point $w=(y, t)\in \mathcal{O}_\Gamma (M)$, we denote the $\mathbb{R}$ coordinate $t$ of $w$ by $t(w)$.
\par
There exist unique group elements $\gamma_{z_C}, \gamma_z, \gamma_x\in \Gamma$ such that
$$z_C = \gamma_{z_C} \cdot P_0^{t(z_C)},$$
$$x = \gamma_{x} \cdot P_0^{t(x)},$$
$$z = \gamma_{z} \cdot P_0^{t(z)}.$$
Hence,
$$x = \gamma_{x} \cdot P_0^{t(x)} =g \cdot \gamma_{z} \cdot P_0^{t(z)}$$
so that $\gamma_x=g \gamma_z$, namely,  $g=\gamma_x \gamma_z^{-1}$ or $g^{-1}=\gamma_z\gamma_x^{-1}$.
\par
Define a trivialization on $C$:
$$t_C : (E\oplus_p H)_{x\in C} \longrightarrow C \times (E\oplus_p H) $$
in such a way that, for each $x\in C$, an affine isometry is given by the formula:
\[
\begin{array}{lccc}
t_C(x): & (E\oplus_p H)_x & \longrightarrow & (E\oplus_p H)_{z_C}= E\oplus_p H \\
 & (\xi, \eta) & \longmapsto & \Big( \gamma_{z_C}^{-1} \gamma_z \xi, \;\;  \gamma_z\gamma_x^{-1} \eta \Big),
\end{array}
\]
where, by abusing notation, the first component $\gamma_{z_C}^{-1} \gamma_z \xi$ represent the affine isometric action $\alpha$ of the group element $\gamma_{z_C}^{-1} \gamma_z$ on $\xi\in E$, and the second component $\gamma_z\gamma_x^{-1} \eta$ represent the image of a vector $\eta\in H$ under the unitary operator $U_{\gamma_z\gamma_x^{-1}}$ on the Hilbert space $H$ given by the linearization of $(\Gamma, M)$.
\par
Finally, we check the requirements (1) and (2) in Definition \ref{fibred} of fibred coarse embeddability for the above constructioins.
\par
{\it The requirement (1) in Definition \ref{fibred}.}  For any $x, x'\in C$, there exist unique $z'\in \text{Ball}_{\mathcal{O} M} (z_C, 3R)$ and $g'\in \Gamma$ such that $g'z'=x'$. Since $\Gamma$ acts on $\mathcal{O} M$ by restrictions of unitary operators on the Hilbert space $H$, we have
$$\delta_\Gamma (x, x') = |g g'^{-1}|+ \|z-z'\|_H.$$
Consequently, we have
\[
\begin{array}{rl}
   & \Big\| t_C(x) (s(x)) -t_C (x') (s (x')) \Big\|_p  \\
 = & \Big\| t_C(x) ( (b_{\gamma_x^{-1}}, x) ) - t_C(x') ( (b_{\gamma_{x'}^{-1}}, x') ) \Big\|_p  \\
 = & \Big( \| \gamma_{z_C}^{-1} \gamma_z b_{\gamma_x^{-1}} - \gamma_{z_C}^{-1} \gamma_{z'} b_{\gamma_{x'}^{-1}}\|_E^p + \|\gamma_z\gamma_x^{-1}x - \gamma_{z'}\gamma_{x'}^{-1}x'\|_H^p \Big)^{\frac{1}{p}} \\
 = & \Big(  \|b_{\gamma_z \gamma_x^{-1}}- b_{\gamma_{z'}\gamma_{x'}^{-1}} \|_E^p +   \|g^{-1}x - g'^{-1} x' \|_H^p        \Big)^{\frac{1}{p}}  \\
 = & \Big( \|b_{g^{-1}}-b_{g'^{-1}}\|_E^p + \|z-z'\|_H^p \Big)^{\frac{1}{p}} \\
 \in & \big[ \rho_1, \rho_2 \big] \Big(  |gg'^{-1}| + \|z-z'\|_H \Big) \\
 = & \big[ \rho_1, \rho_2 \big] \Big( \delta_\Gamma (x, x') \Big)
\end{array}
\]
by the formula ($*$) above. That is, the requirement (1) is fulfilled.
\par
{\it The requirement (2) in Definition \ref{fibred}.} Suppose we have subsets $C_1, C_2\subset \mathcal{O}_\Gamma (M)_{[t_R, \infty)}$ of diameter less than $R$ such that $C_1\cap C_2\not= \emptyset$. Let $x\in C_1\cap C_2$.
Then the affine isometries at $x$ in trivializatioins over $C_1$ and $C_2$ are respectively the following mappings:
\[
\begin{array}{lccc}
t_{C_1}(x): & (E\oplus_p H)_x & \longrightarrow & (E\oplus_p H)_{z_{C_1}}= E\oplus_p H \\
 & (\xi, \eta) & \longmapsto & (\gamma_{z_{C_1}}^{-1} \gamma_{z_1} \xi, \;\; g_1^{-1} \eta), \\
t_{C_2}(x): & (E\oplus_p H)_x & \longrightarrow & (E\oplus_p H)_{z_{C_2}}= E\oplus_p H \\
 & (\xi, \eta) & \longmapsto & (\gamma_{z_{C_2}}^{-1} \gamma_{z_2} \xi, \;\; g_2^{-1} \eta).
\end{array}
\]
Note that there exists $\gamma\in \Gamma$ such that $\gamma z_1 =z_2$. Since $g_1 z_1 =x = g_2z_2$, we have $\gamma=g_2^{-1}g_1$. On the other hand, $z_1 = \gamma_{z_1} P_0$, $z_2 = \gamma_{z_2} P_0$. It follows that
$\gamma\gamma_{z_1}P_0=\gamma_{z_2}P_0$, or equivalently, $\gamma^{-1}=\gamma_{z_1} \gamma_{z_2}^{-1}$. Consequently,
\[
\begin{array}{lll}
&   &  \big[ t_{C_1}(x) \circ ( t_{C_2}(x) ) ^{-1} \big] \big( (\xi, \eta) \big)   \\
& = &  \big( (\gamma_{z_{C_1}}^{-1} \gamma_{z_1} ) (\gamma_{z_{C_2}}^{-1} \gamma_{z_2} )^{-1} \xi, \;\; g_1^{-1} g_2 \eta \big) \\
& = &  \big( (\gamma_{z_{C_1}}^{-1}  \gamma^{-1}  \gamma_{z_{C_2}} \xi, \;\; \gamma^{-1} \eta  )  \big)
\end{array}
\]
Hence, we arrive at an affine isometry for $C_1\cap C_2$:
\[
\begin{array}{lccc}
t_{C_1C_2}: & E\oplus_p H & \longrightarrow & E\oplus_p H \\
 & (\xi, \eta) & \longmapsto & \Big( \gamma_{z_{C_1}}^{-1} \gamma^{-1} \gamma_{z_{C_2}} \xi,  \;\;\; \gamma^{-1} \eta \Big),
\end{array}
\]
which is independent of particular elements $x\in C_1\cap C_2$, such that $t_{C_1}(x)\circ t_{C_2}^{-1}(x)=t_{C_1C_2}$ for all $x\in C_1\cap C_2$. That is, the requirement (2) is fulfilled. This completes the whole proof.
\qed

\section{From coarse geometry of warped cones to proper group actions}
In this section, we seek certain converse implications of the previous sections, i.e., fibred coarse embeddability of warped cones implies proper affine isometric actions of the discrete groups on
various Banach spaces. The warped cones in this section are defined as in Definition \ref{warped cones on general}.
\par
The main result of this section is as follows:
\par
\begin{theorem}
\label{cones to groups}
Let $M$ be a compact metric space and let $\Gamma$ be a finitely generated group acting on $M$ freely by isometries. Assume $M$ is endowed with a finite $\Gamma$-invariant regular Borel measure $\mu$. If the warped cone $\mathcal{O}_\Gamma (M)$ admits a fibred coarse embedding in a Banach space $E$, then $\Gamma$ admits a proper affine isometric action on the ultraproduct $\big( L_p(M, \mu; E) \big)_{\mathcal{U}}$, for any $p\in [1, +\infty)$ and any non-principal ultrafilter $\mathcal{U}$ on the natural numbers $\mathbb{N}$, where $L_p(M, \mu; E)$ is the Lebesgue-Bochner spaces.
\end{theorem}
\par
The proof is motivated by a recent work of S. Arnt \cite{Arnt}, where in particular it is observed that $R$-local proper affine isometric actions of the group pass to a proper affine isometric actions of the whole group
by taking ultralimits. This is made precise as follows.
\par
Let $\Gamma$ be a finitely generated group, and $R>0$ be a real number.
\par
Let $X$ be a set. A map $\alpha: \Gamma\times X \to X$ is said to be a $R$-local action of $\Gamma$ on $X$ if: (1) for all $g\in \Gamma$ such that
$d(e, g)<R$, $\alpha(g): X\to X$ is a bijection; (2) for all $g, h\in \Gamma$ such that $d(e, g)$, $d(e, h)$ and $d(e, gh)$ are less that $R$, $\alpha(gh)=\alpha(g)\alpha(h)$.
\par
Let $E$ be a Banach space. A map $\pi: \Gamma\times E\to E$ is said to be a $R$-local isometric representation of $\Gamma$ on $E$, if $\pi$ is a $R$-local action of $\Gamma$ on $E$ and for all $g\in \Gamma$ such that $d(e, g)<R$,
$\pi(g): E\to E$ is a linear isometry. In this case, a map $b: \Gamma\to E$ such that, for all $g, h\in \Gamma$ such that  $d(e, g)$, $d(e, h)$ and $d(e, gh)$ are less that $R$, $b(gh)=\pi(g)b(h)+b(g)$, is called a $R$-local cocycle with respect to $\pi$.
\par
Let $E$ be a Banach space. A map $\alpha: \Gamma\times E \to E$ is said to be a $R$-local affine isometric action of $\Gamma$ on $E$ if it can be written as $\alpha(g) (\cdot) =\pi(g) (\cdot) + b(g)$ where $\pi$ is a $R$-local isometric representation and $b$ is a $R$-local cocycle with respect to $\pi$.
\par
\begin{lemma} [S. Arnt \cite{Arnt}]
\label{takingultralimit}
Let $\Gamma$ be a finitely generated group and $(E_R)_{R\in \mathbb{N}}$ be a family of Banach spaces. Let $E_{\mathcal{U}}$ be the ultraproduct of the family $(E_R)_{R\in \mathbb{N}}$ with respect to a non-principal ultrafilter
$\mathcal{U}$ on $\mathbb{N}$.  Assume that, for each $R\in \mathbb{N}$, $\Gamma$ admits a $R$-local affine isometric action $\alpha_R$ on $E_R$ with
$\alpha_R(g)(\cdot)=\pi_R(\cdot)+b_R$. If for all $g\in \Gamma$, the sequence $(b_R(g))_{R\in \mathbb{N}}$ gives rise to an element in $E_{\mathcal{U}}$, then, by passing to ultralimit, there exists an affine isometric action $\alpha$ of the whole group $\Gamma$ on the
ultraproduct $E_{\mathcal{U}}$.
\end{lemma}
\par
Now we are ready to prove the main theorem of this section.
\par
\par
\noindent {\bf Proof of Theorem \ref{cones to groups}.} Without loss of generality, assume $\mu (M)=1$. For $t>1$, denote by $\mathcal{O}_\Gamma (M)_t$ the $t$-level set of the warped cone $\mathcal{O}_\Gamma (M)$. This is the space $M$ with the warped metric $\delta_\Gamma$ which is obtained by warping the metric $td_M$ along the $\Gamma$-actions.
\par
For any $R>0$, denote $\Gamma_R=\{\gamma\in \Gamma: |\gamma|<R\}$. Since $M$ is compact and the $\Gamma$-action is free, there exists $t_R>0$ large enough such that $\delta_\Gamma (x, \gamma x)=|\gamma|$ for all $x\in \mathcal{O}_\Gamma (M)_{t_R}$ and $\gamma\in \Gamma_R$. Since the $\Gamma$-action is by isometries,  one can further choose $t_R$ large enough such that, for any $x\in \mathcal{O}_\Gamma (M)_{t_R}$,
$$\text{Ball}_{\mathcal{O}_{\Gamma}(M)} \big( x, 10R \big) \;  \subset \bigsqcup _{|\gamma|\leq 10R} {\text{Ball}_{\mathcal{O} M}} \big(\gamma x, 10R \big).$$
\par
Since the warped cone $\mathcal{O}_\Gamma (M)$ admits a fibred coarse embedding into the Banach space $E$ as in Definition \ref{fibred}, one can choose $t_R$ even larger such that for any subset $C\subset \mathcal{O}_\Gamma (M)_{t_R}$ of diameter less than $R$, there exists a trivialization
$$ t_{N_{10R}(C)} : (E_x)_{x\in N_{10R}(C)} \longrightarrow N_{10R}(C)\times E $$
consisting of affine isometries  $ t_{N_{10R}(C)} (x) : E_x \to E$ for all $x\in N_{10R}(C)$ satisfying the conditions (1) and (2) in Definition \ref{fibred}, where
$$N_{10R}(C):=\{x\in \mathcal{O}_\Gamma (M) \; | \; \delta_\Gamma (x, C) < 10 R \} $$
denotes the $10R$-neighborhood of $C$ in the warped cone $ \mathcal{O}_\Gamma (M) $.
\par
Take a {\em finite} Borel partition $\mathcal{C}_R$ of $\mathcal{O}_\Gamma (M)_{t_R}$ such that, for any $C\in \mathcal{C}_R$,
$$\text{diameter}(C)<R$$
in $\mathcal{O}_\Gamma (M)_{t_R}$. For any $C', C\in \mathcal{C}_R$ with $N_{10R}(C') \cap N_{10R}(C) \not= \emptyset$, there exists an affine isometry
$$t_{N_{10R}(C'), N_{10R}(C)}: E \to E$$
as appeared in the condition (2) in Definition \ref{fibred}. Denote by $u_{C', C}$ the linear part of $t_{N_{10R}(C'), N_{10R}(C)}$. Then, for any vectors $v, w\in E$,
$$\Big( u_{C', C} \Big) (v-w) =  \Big(  t_{N_{10R}(C'), N_{10R}(C)} \Big)  (v) -  \Big(  t_{N_{10R}(C'), N_{10R}(C)} \Big)  (w)              , $$
and for any $C'', C', C\in \mathcal{C}_R$ with
$N_{10R}(C'') \cap N_{10R}(C') \cap N_{10R}(C) \not= \emptyset$, we have
$$u_{C'', C} = u_{C'', C'} \; u_{C', C}.$$
\par
In considering the Lebesgue-Bochner space $L_p(M, \mu; E)$, we identify $\mathcal{O}_\Gamma (M)_{t_R}$ with $M$, and denote by $\overline{u}_{C', C}: L_p(M, \mu; E) \to L_p(M, \mu; E)$ the linear isometry induced by $u_{C', C}$:
$$ \Big( \overline{u}_{C', C} (\xi) \Big) (x) = u_{C', C} \Big( \xi (x) \Big),$$
for all $\xi\in L_p(M, \mu; E)$ and $x\in M$. Moreover, for any $C\in \mathcal{C}_R$ and $\xi\in L_p(M, \mu; E)$, we denote by $\xi_{|C}$ the restriction of $\xi$ on $C$, i.e.
\[
\xi_{|C} = \left\{
\begin{array}{rl}
\xi(x), & x\in C; \\
0, & x\in M\backslash C.
\end{array}
\right.
\]
\par
Now, we define a $R$-local affine isometric action $\alpha^R$ of $\Gamma_R$ on $L_p(M, \mu; E)$ by the formula
$$\alpha^R_\gamma (\xi) = \pi ^R_\gamma (\xi) + b^R_\gamma,$$
for all $\gamma\in \Gamma_R$ and $\xi\in L_p(M, \mu; E)$, in which
\[
\begin{array}{rcl}
\pi^R_\gamma(\xi)   & = &  \bigoplus_{C', C\in \mathcal{C}_R}  \overline{u}_{C', C} \big( (\lambda_\gamma (\xi_{|C}) )_{|C'} \big), \\
b^R_\gamma   &  =  &  \bigoplus_{C\in \mathcal{C}_R} (b^R_\gamma)_{|C},
\end{array}
\]
where $(b^R_\gamma)_{|C} \in L_p(C, \mu; E)$ is defined by
$$ \Big( (b^R_\gamma)_{|C}\Big)  (x)= \Big( t_{N_{10R}(C)}(x) \Big) (s(x)) -  \Big( t_{N_{10R}(C)}(\gamma^{-1}x) \Big) (s(\gamma^{-1}x))$$
for all $x\in C$ for any $C\in \mathcal{C}_R$,  and $\gamma\in \Gamma_R$.
\par
In the rest of the proof, we verify that these $\alpha^R$ are indeed $R$-local affine isometric actions as required. To begin with, it is easy to see that  $\pi^R_\gamma$ is a linear isometry for any $\gamma\in \Gamma_R$. Indeed,
\[
\begin{array}{lll}
\| \pi^R_\gamma (\xi) \|^p_p & = & \sum_{C'\in \mathcal{C}_R}\sum_{C\in \mathcal{C}_R} \| \overline{u}_{C', C} ( ( \lambda_\gamma (\xi_{|C}) )_{|C'} ) \|^p_p  \\
& = & \sum_{C\in \mathcal{C}_R}  \| \lambda_\gamma (\xi_{|C}) \|^p_p   \\
& = & \sum_{C\in \mathcal{C}_R}  \| \xi_{|C} \|^p_p   \\
& = & \|\xi\|^p_p.
\end{array}
\]
Similarly, since for all $\gamma\in \Gamma_R$ and all $x\in C$ for any $C\in \mathcal{C}_R$, we have
$$\rho_1 (|\gamma|) \leq \big\| \big( (b^R_\gamma)_{|C}\big)  (x) \big\|_E \leq \rho_2 (|\gamma|), $$
it follows that, for all $\gamma\in \Gamma_R$,
$$\rho_1 (|\gamma|) \leq \|b^R_\gamma \|_{L_p(M, \mu; E)} \leq \rho_2 (|\gamma|).   \eqno{(**)} $$
\par
{\bf Claim 1.} $\pi^R$ is a $R$-local linear isometric representation of $\Gamma$ on $L_p(M, \mu; E)$, i.e. $\pi^R_{\gamma_2\gamma_1}=\pi^R_{\gamma_2}\pi^R_{\gamma_1}$ for all $\gamma_2, \gamma_1\in \Gamma_R$ with $\gamma_2\gamma_1\in \Gamma_R$.
\par
Indeed, for any $\gamma_2, \gamma_1\in \Gamma_R$ with $\gamma_2\gamma_1\in \Gamma_R$ and $\xi \in L_p(M, \mu; E)$, we have
\[
\begin{array}{rcl}
\pi^R_{\gamma_2}\big( \pi^R_{\gamma_1}(\xi)\big) & = & \bigoplus_{C', C\in \mathcal{C}_R } \pi^R_{\gamma_2} \Big( \overline{u}_{C', C} \big( \big( \lambda_{\gamma_1} (\xi_{|C}) \big)_{|C'} \big) \Big) \\
 & = & \bigoplus_{C', C\in \mathcal{C}_R } \Big\{ \bigoplus_{C''\in \mathcal{C}_R} \overline{u}_{C'', C'}  \Big\{ \lambda_{\gamma_2} \Big[ \big[ \overline{u}_{C', C} \big( \big( \lambda_{\gamma_1} (\xi_{|C})
                \big)_{|C'} \big) \big]_{|C''} \Big] \Big\} \Big\} \\
 & = & \bigoplus_{C'', C', C\in \mathcal{C}_R} \overline{u}_{C'', C'}\; \overline{u}_{C', C} \; \Big[ \lambda_{\gamma_2} \Big[ \Big( \big( \lambda_{\gamma_1} (\xi_{|C}) \big)_{|C'} \Big)_{|C''} \Big] \Big] \\
 & = & \bigoplus_{C'', C\in \mathcal{C}_R} \overline{u}_{C'', C} \; \Big( \big( \lambda_{\gamma_2\gamma_1} (\xi_{|C}) \big)_{|C''} \Big) \\
 & = & \pi^R_{\gamma_2\gamma_1} (\xi).
\end{array}
\]
\par
{\bf Claim 2.} The map $b^R: \Gamma_R\to L_p(M, \mu; E)$ is a $R$-local $\pi^R$-cocycle, i.e.
$$b^{R}_{\gamma_2\gamma_1} =\pi^R_{\gamma_2} (b^R_{\gamma_1}) + b^R_{\gamma_2}$$
for all $\gamma_2, \gamma_1\in \Gamma_R$ with $\gamma_2\gamma_1\in \Gamma_R$.
\par
Indeed, we have
\begin{itemize}
\item  $b^R_{\gamma_2\gamma_1}=\bigoplus_{C'\in \mathcal{C}_R} (b^R_{\gamma_2\gamma_1})_{|C'}$, where for any $x\in C'\in \mathcal{C}_R$,
                $$ \Big( (b^R_{\gamma_2\gamma_1})_{|C'} \Big) (x) = \Big( t_{N_{10R}(C')}(x)\Big)  (s(x)) -  \Big( t_{N_{10R}(C')}({(\gamma_2\gamma_1)}^{-1} x) \Big) (s({(\gamma_2\gamma_1)}^{-1}x));$$
\item  $b^R_{\gamma_2}=\bigoplus_{C'\in \mathcal{C}_R} (b^R_{\gamma_2})_{|C'}$, where for any $x\in C'\in \mathcal{C}_R$,
                $$ \Big( (b^R_{\gamma_2})_{|C'} \Big) (x) = \Big( t_{N_{10R}(C')}(x) \Big)  (s(x)) -  \Big( t_{N_{10R}(C')}({\gamma_2}^{-1} x) \Big) (s({\gamma_2}^{-1}x));$$
\item  $b^R_{\gamma_1}=\bigoplus_{C\in \mathcal{C}_R} (b^R_{\gamma_1})_{|C}$, where for any $x\in C\in \mathcal{C}_R$,
                $$ \Big( (b^R_{\gamma_1})_{|C} \Big) (x) =  \Big( t_{N_{10R}(C)}(x) \Big)  (s(x)) -  \Big(  t_{N_{10R}(C)}({\gamma_1}^{-1} x)  \Big) (s({\gamma_1}^{-1}x)).$$
\end{itemize}
\par
We also have
$$\pi^R_{\gamma_2} (b^R_{\gamma_1}) = \bigoplus_{C', C\in \mathcal{C}_R} \overline{u}_{C', C} \Big( \Big(\lambda_{\gamma_2} \big(  (b^R_{\lambda_1})_{|C} \big)\Big)_{|C'} \Big), $$
so that, for any $C'\in \mathcal{C}_R$,
$$ \Big( \pi^R_{\gamma_2} (b^R_{\gamma_1}) \Big)_{|C'} = \bigoplus_{C\in \mathcal{C}_R} \overline{u}_{C', C} \Big( \Big(\lambda_{\gamma_2} \big(  (b^R_{\lambda_1})_{|C} \big)\Big)_{|C'} \Big).$$
\par
For any $x\in C'$, there exits a unique $C\in \mathcal{C}_R$ such that $\gamma_2^{-1}x\in C$. It follows that, for these $C'$ and $C$, we have
\[
\begin{array}{rl}
    & \Big( \pi^R_{\gamma_2} (b^R_{\gamma_1}) \Big)_{|C'} \big( x \big)  \\
  = & u_{C', C} \Big( b^R_{\gamma_1} (\gamma_2^{-1} x ) \Big)   \\
  = & u_{C', C} \Big( \big( t_{N_{10R}(C)}(\gamma_2^{-1}x) \big) (s(\gamma_2^{-1}x)) -  \big( t_{N_{10R}(C)}({\gamma_1}^{-1} \gamma_2^{-1}x) \big) (s({\gamma_1}^{-1}\gamma_2^{-1}x))  \Big) \\
  = & t_{N_{10R}(C'), N_{10R}(C)}  \Big( \big( t_{N_{10R}(C)}(\gamma_2^{-1}x) \big) (s(\gamma_2^{-1}x)) -  \big( t_{N_{10R}(C)}({\gamma_1}^{-1} \gamma_2^{-1}x) \big) (s({\gamma_1}^{-1}\gamma_2^{-1}x))  \Big) \\
  = & \Big( t_{N_{10R}(C')}(\gamma_2^{-1}x) \Big) (s(\gamma_2^{-1}x)) -  \Big( t_{N_{10R}(C')}({(\gamma_2\gamma_1)}^{-1} x) \Big) (s({(\gamma_2\gamma_1)}^{-1} x)).
\end{array}
\]
It follows that, for all $x\in C'\in \mathcal{C}_R$,
$$\Big( b^{R}_{\gamma_2\gamma_1} \Big)_{|C'} (x)  = \Big( \pi^R_{\gamma_2} (b^R_{\gamma_1}) \Big)_{|C'} (x) + \Big( b^R_{\gamma_2} \Big)_{|C'} (x).$$
This proves Claim 2.
\par
Note that we already have the controlling inequalities $(**)$, i.e., for all $\gamma\in \Gamma_R$,
$$\rho_1 (|\gamma|) \leq \|b^R_\gamma \|_{L_p(M, \mu; E)} \leq \rho_2 (|\gamma|).  $$
Passing to ultralimit, by using Lemma 4.2, we complete the proof.  \qed

\par

Recall that a discrete group $\Gamma$ has Gromov's a-T-menability \cite{CCJJV}, or the Haagerup property, if $\Gamma$ admits a proper affine isometric action on a Hilbert space. The advantage of working with Hilbert space over other $L_p$-spaces is the availability of the effective tools of conditionally negative definite functions or kernels. A group $\Gamma$ is a-T-menable if and only if there exists a conditionally negative definite functions $\psi: \Gamma\to \mathbb{R}^+$ which is metrically proper in the sense that $\lim_{\gamma\to\infty} \psi(\gamma)=\infty$. This allows us to relax the condition of isometric actions in the above Theorem \ref{cones to groups} to any actions by homeomorphisms. With only a minor modifications of the proof in [\cite{Roe05}, Proposition 4.4] (see also [\cite{Saw}, Proposition 6.1] and [\cite{CWW}, Proposition 2.10]), we have the following strengthened result.
We include the proof for convenience of readers.
\par
\begin{proposition}
Let $M$ be a compact metric space and let $\Gamma$ be a finitely generated group acting on $M$ by homeomorphisms. Assume that $M$ is endowed with a finite $\Gamma$-invariant measure $\mu$ and there exits a subset $P\subset M$ of positive measure on which the action of $\Gamma$ is free. If the warped cone $\mathcal{O}_\Gamma (M)$ admits a fibred coarse embedding into Hilbert space, or more generally, $\mathcal{O}_\Gamma (M)$ is asymptotically (coarse) embeddable into Hilbert space in the sense of Rufus Willett in \cite{Will15}, then the group $\Gamma$ has the Haagerup property.
\end{proposition}
\par
{\bf Proof.} Without loss of generality in the proof, assume that $P=M$ and $\mu(M)=1$. Suppose $\mathcal{O}_\Gamma (M)$ admits a fibred coarse embedding into Hilbert space $H$, as in Definition \ref{fibred} with $E=H$. For any $R>0$, there exits $t_R$ large enough such that for all $x\in \mathcal{O}_\Gamma (M)_{t_R}$ and $\gamma\in \Gamma$ with $|\gamma |<R$,
$$\delta_\Gamma (x, \gamma x)=|\gamma|, $$
and for any $C\subset \mathcal{O}_\Gamma (M)_{t_R}$ of diameter less $R$, there exists a trivialization $t_C$ consisting of a family of affine isometries $t_C (x) :H\to H $ for all $x\in C$ satisfying the conditions (1) and (2) in Definition . Hence, we obtain a well-defined $R$-local conditionally negative definite kernel $\mathcal{O}_\Gamma (M)_{t_R}$:
$$K^R(x, y)=\| t_C(x)(s(x)) -t_C(y)(s(y)) \|^2 $$
such that
$$\rho_1(d(x, y)) \leq K^R(x,y) \leq \rho_2 (d(x,y))$$
for all $x, y\in \mathcal{O}_\Gamma (M)_{t_R}$ with $\delta_\Gamma (x, y)\leq R$.
\par
For each $\gamma\in \Gamma$ with $|\gamma |\leq R$, define
$$h^R(\gamma)=\int_M K^R(x, \gamma x) d\mu(x).$$
Then $h^R$ is a $R$-local (right-hand) conditionally negative definite function on $\Gamma$, i.e., the kernel $(\gamma_1, \gamma_2)\mapsto h^R (\gamma_1\gamma_2^{-1})$ is a $R$-local conditionally negative definite functions on $\Gamma\times \Gamma$. Moreover, for all $\gamma\in \Gamma$ with $|\gamma|\leq R$,
$$\rho_1(|\gamma|)\leq h^R(\gamma)\leq \rho_2(|\gamma|).$$
By passing to a subsequence and taking limit, we can arrive at a conditionally negative definite function $h(\gamma)=\lim_{k} h^{R_k}(\gamma)$ on the whole $\Gamma$ which is metrically proper. This is the Haagerup property of $\Gamma$. \qed
\par

\section{Finitely generated dense subgroups of compact metrizable groups}
The results of the previous sections apply very well to the case of finitely generated dense subgroups of compact metrizable groups. Note that the Banach spaces involved in
Theorem \ref{group to cone} and Theorem \ref{cones to groups} can be replaced by Banach spaces in any other class $\mathfrak{B}$ of Banach spaces which is closed under taking $\ell_p$-direct sums,
the Lebesgue-Bochner $L_p$-spaces and ultraproducts, i.e.
$$E, F\in \mathfrak{B} \Longrightarrow E\oplus_p F\in \mathfrak{B}, \;\; L_p(\Omega, \mu; E)\in \mathfrak{B}, \;\; (E)_{\mathcal{U}} \in \mathfrak{B}.$$
For example, the class $\mathfrak{L}_p$ of all $L_p$-spaces , the class of uniformly convex Banach spaces, or the class of Banach spaces with non-trivial type. We have the following summarization:
\par
\begin{theorem}
\label{compactmetrizablegroups}
Let $G$ be a compact metrizable group, and $\Gamma$ be a finitely generated dense subgroup of $G$ acting on $G$ by left translations. Then
\begin{enumerate}
\item[(1)] $\mathcal{O}_\Gamma (G)$ admits a fibred coarse embedding into a $L_p$-space $(1\leq p <\infty)$ if and only if $\Gamma$ admits a proper affine isometric action on a $L_p$-space.
\item[(2)] $\mathcal{O}_\Gamma (G)$ admits a fibred coarse embedding into a uniformly convex Banach space if and only if $\Gamma$ admits a proper affine isometric action on a uniformly convex Banach space.
\item[(3)] $\mathcal{O}_\Gamma (G)$ admits a fibred coarse embedding into a Banach space with non-trivial type if and only if $\Gamma$ admits a proper affine isometric action on a Banach space with non-trivial type.
\item[(4)] $\mathcal{O}_\Gamma (G)$ admits a fibred coarse embedding into Hilbert space if and only if $\Gamma$ has the Haagerup property.
\end{enumerate}
In particular, if $\Gamma$ is a finitely generated dense subgroup of a compact Lie group $G$, then the above statements hold.
\end{theorem}
\par

We discuss now some more concrete examples or counterexamples.
\par

Firstly, there are warped cones over a-T-menable groups that do not embed coarsely into Hilbert space, or more generally, into any $L_p$-space, $1\leq p<\infty$. J. Bourgain and A. Gamburd proved \cite{BG08} that for
appropriately chosen subgroups $\mathbb{F}_n$ in $\text{SU(2)}$, the action of $\mathbb{F}_n$ on $\text{SU(2)}$ has a spectral gap. They further generalized this result to $\text{SU(d)}$ in \cite{BG12}.
P. W. Nowak and D. Sawicki \cite{NowSaw} very recently showed that, due to spectral gaps, the associated warped cones $\mathcal{O}_{\mathbb{F}_n}(\text{SU(d)})$ is not coarsely embeddable into any $L_p$-space, $1\leq p<\infty$.
However, it follows from Theorem \ref{compactmetrizablegroups} above that $\mathcal{O}_{\mathbb{F}_n}(\text{SU(2)})$ admits a {\em fibred} coarse embedding into Hilbert space.
In general, Y. Cornulier proved in \cite{Corn06} that, if $G$ is a compact Lie group, every countable subgroup of $G$ has the Haagerup property if and only if $G$ is locally isometric to $\text{SO}(3)^k \times \text{U}(1)^m$
for some integers $k,m$. Hence, the corresponding warped cones admit a {\em fibred} coarse embedding into a Hilbert space.

\par
Secondly, examples of dense property (T) subgroups of compact Lie groups can be found for instance with $G=\text{SO}(5)$ \cite{Lub}. It was shown early in \cite{Roe05} that $\mathcal{O}_\Gamma ( \text{SO}(5) )$
does not coarse embed into a Hilbert space. It follows from the above Theorem \ref{compactmetrizablegroups} that, even more, $\mathcal{O}_\Gamma ( \text{SO}(5) )$ does not admit a fibred coarse embedding into a
Hilbert space. In general, Y. Cornulier proved in \cite{Corn08} that a compact Lie group admits a countable dense subgroup with property (T) if and only if no quotient of $G$ is isomorphic either to $\text{SO}(3)$
nor to $\text{U}(1)$. The corresponding warped cones do not admit a fibred coarse embedding into Hilbert space.

\par
Thirdly, consider a finitely generated group $\Gamma$ with strong Banach property (T) \cite{Laf}, that is, Property (T) with respect to all Banach spaces of non-trivial type. Examples of such groups include lattices in $\text{SL}_3(\mathbb{Q}_p)$, as shown by V. Lafforgue  \cite{Laf}. This class was extended by B. Liao in \cite{Liao} to connected almost $F$-simple algebraic groups whose $F$-split rank is at least $2$, where $F$ denotes a non-Archimedean local field. Let $Y$ be a compact metric space on which $\Gamma$ acts freely by measure-preserving isometries.
It was shown in \cite{NowSaw} that the warped cone $\mathcal{O}_\Gamma (Y)$ for $\Gamma$ and $Y$ does not embed coarsely
into any Banach space with non-trivial type. It follows from Theorem \ref{cones to groups} that actually $\mathcal{O}_\Gamma (Y)$ does not admit a fibred coarse embedding into any Banach space with non-trivial type.

\par
In the rest of this section, let $\Gamma=\langle S \rangle  $ be a finitely generated residually finite group,  with $\{\Gamma_n\}_{n\in \mathbb{N}}$ a nested decreasing sequence of finite index normal subgroups of $\Gamma$ with trivial intersection. We have the following sequence of epimorphisms
$$1\leftarrow \Gamma/\Gamma_1 \leftarrow \Gamma/\Gamma_2 \leftarrow \cdots.$$
The inverse limit of this sequence, denoted by $\widehat{\Gamma} ((\Gamma_n))$, is call the {\em profinite completion of $\Gamma$ with respect to the sequence $\{\Gamma_n\}$}. This is a subset of the Cartesian product
$\prod_{n\in \mathbb{N}} \Gamma/\Gamma_n$:
$$\widehat{\Gamma} ((\Gamma_n))=\{ (g_n)_{n\in \mathbb{N}} | g_{n-1}= q_n (g_n), \; \forall n\geq 1 \},$$
where $q_n: \Gamma/\Gamma_n\to \Gamma/\Gamma_{n-1}$ are the intermediate quotient maps. It follows that $\widehat{\Gamma} ((\Gamma_n))$ is a compact group, and the quotient maps $\pi_n: \Gamma\to \Gamma/\Gamma_n$ give rise to a homomorphism from $\Gamma$ to $\widehat{\Gamma} ((\Gamma_n))$ with dense image. The kernel of this map is $\cap_n \Gamma_n$. If the intersection is trivial, then the action of $\Gamma$ on $G$ by left multiplication is free. The group $\widehat{\Gamma} ((\Gamma_n))$ is metrizable and admits a $\Gamma$ invariant measure, the Haar measure. A particular metric on $\widehat{\Gamma} ((\Gamma_n))$ is given by:
$$d\big( (g_n), (g'_n) \big) =\max \{ a_n \cdot |g_ng_n'^{-1}|_S \}, $$
where $a_n$ is a sequence of positive numbers decreasing to zero in such a way that $a_{n+1}< a_{n} / \text{diam} (\Gamma/\Gamma_n)$ \cite{Saw}.

\par
\begin{theorem}[Sawicki, \cite{Saw}]
Let $G=\widehat{\Gamma} ((\Gamma_n))$ be the profinite completion of a finitely generated residually finite group $\Gamma$ as above. TFAE:
\begin{itemize}
\item $\Gamma$ is amenable;
\item $\mathcal{O}_\Gamma (G)$ has property A;
\item The box space $\Box_{\Gamma_n} \Gamma $ has property A.
\end{itemize}
\end{theorem}

\par
D. Sawicki \cite{Saw} also proved that coarse embeddability into Hilbert space for the warped cone $\mathcal{O}_\Gamma (G)$ and the box space $\Box_{\Gamma_n} \Gamma $ are equivalent, which in turn implies the Haagerup property of $\Gamma$. Generalizing this, we have the following result concerning fibred coarse embeddability and proper isometric actions.
\par
\begin{theorem}
Let $G=\widehat{\Gamma} ((\Gamma_n))$ be the profinite completion of a finitely generated residually finite group $\Gamma$ as above. TFAE:
\begin{enumerate}
\item[(1)] The group $\Gamma$ admits a proper affine isometric action on a $L_p$-space;
\item[(2)] The warped cone $\mathcal{O}_\Gamma (G)$ fibred coarse embeds into a  $L_p$-space;
\item[(3)] The box space $\Box_{\Gamma_n} \Gamma $ fibred coarse embeds into a  $L_p$-space.
\end{enumerate}
The class of $L_p$-spaces can be replaced by other stable class $\mathfrak{B}$ of Banach spaces.
\end{theorem}
\par
{\bf Proof.} (1) $\Longrightarrow$ (2) by Theorem \ref{group to cone}. (2) $\Longrightarrow$ (3) follows from D. Sawicki \cite{Saw} that the box space $\Box_{\Gamma_n} \Gamma $ coarsely embeds in the warped cone
$\mathcal{O}_\Gamma (G)$. (3) $\Longrightarrow$ (1) follows from S. Arnt's main result in \cite{Arnt}.  \qed
\par

As examples, recall that G. Yu proved in \cite{Yu05} that a Gromov hyperbolic group $\Gamma$ admits a proper affine isometric action on $\ell_p (\Gamma\times\Gamma)$ with $p\geq 2$.
It follows from Theorem \ref{group to cone} that,
if in addition $\Gamma$ is residually finite and $G=\widehat{\Gamma} ((\Gamma_n))$ is as above, than the warped cone $\mathcal{O}_\Gamma (G)$ admits a fibred coarse embeds into $\ell_p (\Gamma\times\Gamma)$.
Note also that B. Nica proved in \cite{Nica} that every non-elementary hyperbolic group $\Gamma$ admits a proper affine isometric action on $L_p(\partial \Gamma \times \partial\Gamma)$, where $\partial \Gamma$ denotes
the Gromov boundary of $\Gamma$ and $p$ is large enough.

\par
For groups with the Haagerup property, we summarize as follows:

\par
\begin{corollary}
Let $G=\widehat{\Gamma} ((\Gamma_n))$ be the profinite completion of a finitely generated residually finite group $\Gamma$ as above. TFAE:
\begin{itemize}
\item The group $\Gamma$ has the Haagerup property;
\item The warped cone $\mathcal{O}_\Gamma (G)$ fibred coarse embeds into Hilbert space;
\item The box space $\Box_{\Gamma_n} \Gamma $ fibred coarse embeds into Hilbert space.
\end{itemize}
\end{corollary}

\section{On higher index problems on warped cones}
In this section, we apply results in previous sections to higher index problems on warped cones. It turns out that the maximal coarse Baum-Connes conjecture and the coarse Novikov conjecture holds for a certain class of warped cones which do not coarsely embed into any $L_p$-space for any $1\leq p<\infty$.
\par
Let $X$ be a proper metric space with bounded geometry. Recall that {\em The (maximal) coarse Baum–Connes conjecture} \cite{Yu00,CWY13} is a geometric analogue of the Baum–Connes conjecture and
provides an algorithm for computing the higher indices of generalized elliptic operators on non-compact spaces. It states that a certain assembly map $\mu$ (resp. $\mu_{max}$) from
$\lim_{d\to \infty} K_*(P_d (X))$ to $K_*(C^*(X))$ or $K_*(C^*_{max} (X))$ is an isomorphism, where $K_*(P_d (X))$ is the locally finite $K$-homology group of the Rips complex
$P_d (X)$, and $K_*(C^*(X))$ and $K_*(C^*_{max}(X))$ are respectively the $K$-theory groups of the Roe algebra and the maximal Roe algebra of $X$. {\em The coarse Novikov conjecture}
states that the reduced assembly map $\mu$ is injective, which provides an algorithm to determine when the higher indices are non-zero. The assembly maps are illustrated in the following commutative diagram
\[
\xymatrix{
  &  K_*(C^*_{\max}(X)) \ar[d]^{\lambda_*}  \\
\displaystyle\lim_{d\to \infty} K_*(P_d(X))  \ar[ur]^{\mu_{\max}}  \ar[r]^{\quad \mu} & K_*(C^*(X))
}
\]
where $\lambda: C^*_{\max}(X)\to C^*(X)$ is the canonical quotient map. These conjectures have many applications in topology and geometry.
\par
\begin{theorem}[\cite{CWY13}]
Let $X$ be a discrete metric space with bounded geometry. If $X$ admits a fibred coarse embedding into Hilbert space, then the maximal coarse Baum-Connes conjecture holds for $X$.
\end{theorem}
\par
We also have the following result which can be proved similarly as in \cite{CWWz,CWY13,CWY15}.
\begin{theorem}
Let $X$ be a discrete metric space with bounded geometry. If $X$ admits a fibred coarse embedding into an $L_p$-space with $1\leq p<\infty$,
then the coarse Novikov conjecture holds for $X$.
\end{theorem}
\par
\begin{corollary}
The maximal coarse Baum-Connes conjecture and the coarse Novikov conjecture hold for the warped cone $\mathcal{O}_{\mathbb{F}_n}(\text{SU(d)})$, which is not coarsely embeddable into any $L_p$-space for any $1\leq p<\infty$.
\end{corollary}

Y. Cornulier proved in \cite{Corn06} that, if $G$ is a compact Lie group, every countable subgroup of $G$ has the Haagerup property if and only if $G$ is locally isometric to $\text{SO}(3)^k \times \text{U}(1)^m$ for some integers
$k,m$. As applications, we have the following results:
\par

\begin{corollary}
Let $G$ be a compact Lie group which is locally isometric to $\text{SO}(3)^k \times \text{U}(1)^m$ for some integers $k,m$. Then for any countable subgroup $\Gamma$ of $G$,
the maximal coarse Baum-Connes conjecture and the coarse Novikov conjecture hold for the warped cone $\mathcal{O}_{\Gamma}(G)$.
\end{corollary}

The coarse Novikov conjecture holds for somewhat wider class of warped cones.

\par
\begin{corollary}
Let $\Gamma$ be a finitely generated hyperbolic group and $G$ be a compact metrizable group containing $\Gamma$ as a dense subgroup. Suppose the warped cone $\mathcal{O}_\Gamma (G)$ has bounded geometry. Then the coarse Novikov conjecture holds for $\mathcal{O}_\Gamma (G)$.
\end{corollary}

In particular, we have the following:
\par
\begin{corollary}
Let $\Gamma$ be a finitely generated residually finite hyperbolic group, and $G=\widehat{\Gamma} ((\Gamma_n))$ the corresponding profinite completion. Then the coarse Novikov conjecture holds for $\mathcal{O}_\Gamma (G)$.
\end{corollary}

\vskip 1cm

\mbox{} \\
Research Center for Operator Algebras, \\
Department of Mathematics, \\
East China Normal University, \\
Shanghai \quad 200241, \\
P. R. China.\\
\\
% E-mail: qwang@math.ecnu.edu.cn


\begin{thebibliography}{99}
%
%
\bibitem{Arnt}
        S. Arnt. Fibred coarse embeddability of box spaces and proper isometric affine actions on $L^p$ spaces. arXiv:1508.05033v1[math.GR] 20 Aug 2015. Bull. Belg. Math. Soc. Simon Stevin Volume 23, Number 1 (2016), 21-32.

\bibitem{BaayenGroot}
        P. C. Baayen and J. De Groot. Linearization of locally compact transformation groups in Hilbert spaces, Math. Systems Theory 2 (1968), 363-379.

\bibitem{BFGM}
        U. Bader; A. Furman; T. Gelander; N. Monod. Property (T) and rigidity for actions on Banach spaces. Acta Math. 198 (2007), no. 1, 57-105.

\bibitem{BG08}
        J. Bourgain and A. Gamburd. On the spectral gap for finitely-generated subgoups of $SU(2)$. Invent. Math. 171(1)83-121, 2008.

\bibitem{BG12}
        J. Bourgain and A. Gamburd. A spectral gap theorem in $\text{SU}(d)$, J. Eur. Math. Soc. (JEMS) 14 (2012), no. 5, 1455-1511, DOI 10.4171/JEMS/337. MR2966656.

\bibitem{CWW}
        X. Chen, Q. Wang, and X. Wang. Characterization of the Haagerup property by fibred coarse embedding into Hilbert space. Bulletin of the London Mathematical Society,
        45(5): 1091-1099, 2013.


\bibitem{CWWz}
        X. Chen, Q. Wang, and Z. Wang. Fibred coarse embedding into non-positively curved manifolds and higher index problem. J. Funct. Anal. 267 (2014), no. 11, 4029-4065.


\bibitem{CWY13}
        X. Chen, Q. Wang, and G. Yu. The maximal coarse Baum-Connes conjecture for spaces which admit a fibred coarse embedding into Hilbert space. Adv. Math. 249 (2013), 88-130.


\bibitem{CWY15}
        X. Chen, Q. Wang, and G. Yu. The coarse Novikov conjecture and Banach spaces with Property (H). J. Funct. Anal. 268 (2015), no. 9, 2754-2786.



\bibitem{CCJJV}
        P. -A. Cherix, M. Cowling, P. Jolissant, P. Julg, A. Valette. Groups with the Haagerup property (Gromov's a-T-menability). Birkh\"{a}user. Progress in Math. 197, 2001.

\bibitem{Corn06}
        Yves de Cornulier. Kazhdan and Haagerup Properties in algebraic groups over local fields. J. Lie Theory 16, 67-82, 2006.

\bibitem{Corn08}
        Yves de Cornulier. Dense subgroups with Property (T) in Lie groups. Comment. Math. Helv. 83(1), 55-65, 2008.

\bibitem{DruNow}
        C. Dru\c{t}u and P. W. Nowak. Kazhdan projections, random walks and ergodic theorems. arXiv:1501.03473v1, math.GR, 2015. To appear in Crelle's Journal.

\bibitem{Grom93}
        M. Gromov.  Asymptotic invariants for infinite groups.  Proc. 1991 Sussex Conference on Geometry Group Theory, LMS Lecture Note Ser.182,
        papge 1--295. Cambridge University Press, 1993.


\bibitem{Grom00}
    	M. Gromov. Spaces and questions. Geom. Funct. Anal., Special Volume, Part I (2000) 118--161. GAFA 2000 (Tel Aviv, 1999)


\bibitem{Grom03}
	   M. Gromov. Random walks in random groups. Geom. Funct. Anal. 13(1) (2003) 73--146.


\bibitem{Hig99}
		N. Higson. Counterexamples to the coarse Baum-Connes conjecture. Unpublished paper, 1999, available on the author's website.

\bibitem{Haag}
        U. Haagerup. An example of a nonnuclear $C^*$-algebra, which has the metric approximation property. Invent. Math. 50(3)(1979) 279--293.



\bibitem{HLS}
		N. Higson, V. Lafforgue, and G. Skandalis.  Counterexamples to the Baum-Connes conjecture. Geom. Funct. Anal. 12(2)  (2002)  330--354.


\bibitem{HK01}
        N. Higson and G. Kasparov. $E$-theory and $KK$-theory for groups which act properly and isometrically on Hilbert space. Invent. Math. 144(1) (2001) 23--74.

\bibitem{KS}
        G. Kasparov and G. Skandalis. Groups acting properly on ``bolic'' spaces and the Novikov conjecture. Annals of Mathematics, 158 (2003) 165--206.

\bibitem{Khu}
        A. Khukhro. Box spaces, group extensions and coarse embeddings into Hilbert space. J. Funct. Anal. 263(1)(2012) 115--128.


\bibitem{Kim}
        H. J. Kim. Coarse equivalences between warped cones. Geom. Dedicata 120 (2006), 19-35.


\bibitem{Laf}
    	V. Lafforgue.  Un renforcement de la propri\'{e}t\'{e} (T). Duke Math. J. 143(3) (2008) 559--602.

\bibitem{Liao}
        B. Liao. Strong Banach property (T) for simple algebraic groups of higher rank. J. Topol. Anal. 6 (2014), no. 1, 75-105.

\bibitem{Lub}
        A. Lubotzky.  Discrete groups, expanding graphs and invariant measures.  Progress in Mathematics, Vol. 125, Birkhauser Verlag, 1994.

\bibitem{Marg}
		G. Margulis. Explicit constructions of expanders. Probl. Peredaci Inf.  9(4) (1973) 71--80.

\bibitem{MooreSchlafly}
        John Douglas Moore and Roger Schlafly. On equivariant isometric embeddings. Math. Z. 173 (1980), no. 2, 119-133.

\bibitem{Mostow1957}
        G. D. Mostow. Equivariant embeddings in Euclidean spaces. Ann. of Math. 65 (1957), 432-446.

\bibitem{Nash}
        John Nash. The imbedding problem for Riemannian manifolds. Ann. of Math. (2) 63 (1956), 20-63.

\bibitem{Nica}
        B. Nica. Proper isometric actions of hyperbolic groups on $L_p$-spaces. Compos. Math. 149 (2013), no. 5, 773-792.

\bibitem{Now05}
        P. W. Nowak. Coarse embeddings of metric spaces into Banach spaces. Proc. Amer. Math. Soc. 133 (2005), no. 9, 2589-2596

\bibitem{Now06}
        P. W. Nowak. On coarse embeddability into lp-spaces and a conjecture of Dranishnikov. Fund. Math. 189 (2006), no. 2, 111-116.

\bibitem{Now06-2}
        P. W. Nowak. Group actions on Banach spaces and a geometric characterization of a-T-menability. Topology Appl. 153 (2006), no. 18, 3409-3412.

\bibitem{Now07}
        P. W. Nowak. Coarsely embeddable metric spaces without Property A. J. Funct. Anal. 252 (2007), no. 1, 126-136.

\bibitem{Now15}
        P. W. Nowak. Group actions on Banach spaces. Handbook of group actions. Vol. II, 121-149, Adv. Lect. Math. (ALM), 32, Int. Press, Somerville, MA, 2015.

\bibitem{NowSaw}
        P. W. Nowak and D. Sawicki. Warped cones and spectral gaps. arXiv:1509.04921v1 [math.MG] 16 Sep 2015. Proceedings of the American Mathematical Society 145 (2017), no. 2, 817-823


\bibitem{NowYu}
        P. W. Nowak and G. Yu. Large scale geometry. EMS Textbooks in Mathematics, European Mathematical Society (EMS), Z\"{u}rich, 2012.

\bibitem{Roe95}
        J. Roe. From foliations to coarse geometry and back. Analysis and geometry in foliated manifolds (Santiago de Compostela, 1994), 195-205, World Sci. Publ., River Edge, NJ, 1995.


\bibitem{Roe03}
        J. Roe. Lectures on Coarse Geometry, volume 31 of University Lecture Series. American Mathematical Society, 2003.

\bibitem{Roe05}
        J. Roe. Warped cones and property A. Geom. Topol. 9 (2005), 163-178.

\bibitem{Saw}
        D. Sawicki. Warped cones over profinite completions. arXiv:1509.04669v1 [math.MG] 15 Sep 2015. J. Topol. Anal., posted on 
        2017, DOI 10.1142/S179352531850019X.
        
\bibitem{Saw17}
        D. Sawicki. Super-expanders and warped cones. arXiv:1704.03865v1 [math.MG] 12 Apr 2017.       

\bibitem{Vigolo}
        Federico Vigolo. Mearsure expanding actions, expanders and warped cones. arXiv:1610.05837v1 [math.GT] 19 Oct 2016.


\bibitem{Will09}
        R. Willett. Some notes on property A. Limits of graphs in group theory and computer science, 191-281, EPFL Press, Lausanne, 2009.

\bibitem{Will15}
        R. Willett. Random graphs, weak coarse embeddings, and higher index theory. J. Topol. Anal., 7(3):361-388, 2015.

\bibitem{WiY1}
		R. Willett and G. Yu.   Higher index theory for certain expanders and Gromov monster groups, I.   Adv. Math.  229(3) (2012) 1380--1416.


\bibitem{WiY2}
		R. Willett and G. Yu.   Higher index theory for certain expanders and Gromov monster groups, II.  Adv. Math.  229(3) (2012) 1762--1803.


\bibitem{Yu00}
        G. Yu.  The coarse Baum-Connes conjecture for spaces which admit a uniform embedding into Hilbert space.  Invent. Math. 139 (2000) 201--240.


\bibitem{Yu05}
        G. Yu.  Hyperbolic groups admit proper affine isometric actions on $\ell^p$-spaces. Geom. Funct. Anal. 15(5)(2005) 1144-1151.



\end{thebibliography}
\end{document}